%% file: DEA_GL_MS.tex
\def \maintitle {Joint Variable Selection for Data Envelope Analysis via Group Sparsity}

\def \subject {subject}

\documentclass[mnsc,nonblindrev]{informs3} 

\OneAndAHalfSpacedXI

\usepackage[authoryear, round, longnamesfirst]{natbib}
\bibpunct[, ]{(}{)}{,}{a}{}{,}%
\def\bibfont{\small}%
\def\bibsep{\smallskipamount}%
\def\bibhang{24pt}%
\usepackage{cite}
\def \keywordlist {Data envelopment analysis (DEA), Group Lasso, Lasso, Group variable selection, Alternating direction method of multipliers (ADMM)}

\TheoremsNumberedThrough     
\ECRepeatTheorems

\EquationsNumberedThrough    

\MANUSCRIPTNO{}

\usepackage{geometry}                
\usepackage{graphicx}
\usepackage{amssymb}
\usepackage{amsmath}
\usepackage{amsfonts}
\usepackage{multirow}
\usepackage{epstopdf}
\usepackage{algorithmic}
\usepackage{algorithm}

\input{rax_pdfsetup}
\include{tony_def}

\newcommand{\ieor}{Department of Industrial Engineering \& Operations Research, Columbia University}

\begin{document}



\RUNTITLE{Joint Variable Selection for DEA}
\TITLE{Joint Variable Selection for Data Envelopment Analysis via Group Sparsity}

\ARTICLEAUTHORS{%
\AUTHOR{Zhiwei (Tony) Qin}
\AFF{\ieor, 120th Street, New York, NY 10027, \EMAIL{tonyqin@gmail.com}} 
\AUTHOR{Irene Song}
\AFF{\ieor, 120th Street, New York, NY 10027, \EMAIL{is2306@columbia.edu}}
} 

\ABSTRACT{%
This study develops a data-driven group variable selection method for data envelopment analysis (DEA), a non-parametric linear programming approach to the estimation of production frontiers.  The proposed method extends the group Lasso (least absolute shrinkage and selection operator) designed for variable selection on (often predefined) groups of variables in linear regression models to DEA models.  In particular, a special constrained version of the group Lasso with the loss function suited for variable selection in DEA models is derived and solved by a new tailored algorithm based on the alternating direction method of multipliers (ADMM).  This study further conducts a thorough evaluation of the proposed method against two widely used variable selection methods -- the efficiency contribution measure (ECM) method and the regression-based (RB) test -- in DEA via Monte Carlo simulations.  The simulation results show that our method provides more favorable performance compared with its benchmarks.
}%

\maketitle

%

\section{Introduction}
\label{sec:intro}
Since its introduction in the seminal paper by Charnes, Cooper and Rhodes \citeyearpar{CCR}, data envelopment analysis (DEA) -- a non-parametric linear programming approach to frontier analysis -- has grown into a powerful quantitative analytical tool in management science and operations research \citep{DEAText2}.  A single comprehensive measure of productive efficiency generated by the method has been popularly served as an essential guidepost for making managerial decisions in practice.  In recent years, we have seen many successful applications of DEA in performance evaluation of economic entities (also known as decision making units (DMUs)) in various fields ranging from non-profit sectors such as hospitals to for-profit sectors such as banks.\footnote{To name a few, by means of DEA, Kuntz and Vera \citeyearpar{KuntsVera} conducted performance analysis of hospitals and Song, Kachani and Kani \citeyearpar{SongKachaniKani} estimated relative efficiency of the firms in the U. S. Information Technology sector and investigated its link to stock performance.  For more applications of DEA, interested readers can refer to \emph{Data Envelopment Analysis and Its Applications to Management} \citep{DEAApplications}.}  Along with its rising popularity, DEA has certainly developed into a widely accepted field of research in its own.

Despite the large number of papers published on DEA,\footnote{According to the literature survey by Liu et al. \citeyearpar{LiuEtAl} the DEA field has accumulated over 4,500 papers in the ISI Web of Science database in the last three decades.} surprisingly scant attention has been paid to variable selection in the literature.  
Variable selection approaches in DEA are often based on experts' opinions, past experience or economic theories, as a matter of fact.  The major concern about these approaches is that they are prone to include irrelevant variables and/or omit relevant variables.  They can, therefore, result in model misspecification.  Several studies have demonstrated significant negative impact that model misspecification has on the accuracy of DEA efficiency estimates.  For instance, Sexton et al. \citeyearpar{SextonEtAl} has investigated the effect of including an irrelevant variable in a DEA model and reported that any variable included in the analysis, in fact, can change the shape and position of the production frontier, which in turn alters the ranking of efficiency estimates.  Similarly, Smith \citeyearpar{Smith} has documented the danger of model misspecification when a relevant variable is omitted from a DEA model.  In addition to model misspecification, attention to variable selection is imperative also because the greater the number of variables included in a DEA model, the higher the dimensionality of the production space and the less discerning the analysis \citep{JenkinsAnderson}.  More specifically, an increase in the number of variables included in a DEA model tends to shift the compared DMUs towards the efficient frontier, thus leading to a decline in DEA's discriminatory power \citep{GolanyRoll, FriedEtAl}.  It is therefore essential to limit the number of variables included in the analysis.  Still, there is no consensus on how best to do this.

In this paper, we propose a data-driven joint variable selection method for DEA.  In particular, we extend the group Lasso (least absolute shrinkage and selection operator) designed for variable selection on (often predefined) groups of variables in linear regression models to DEA models.  We derive a special constrained version of the group Lasso with the loss function suited for variable selection in DEA models and solve it by a new tailored algorithm based on the alternating direction method of multipliers (ADMM).  We conduct a thorough performance evaluation of the proposed method against two of the most widely used variable selection approaches in the DEA literature -- the efficiency contribution measure (ECM) method and the regression-based (RB) test -- by means of Monte Carlo simulations.

\section{Data Envelopment Analysis (DEA)}
\label{sec:DEA}
DEA is a non-parametric mathematical programming approach to the estimation of production frontiers.  It solves a set of linear programs (LPs) to construct a piecewise linear efficient production frontier, which represents the relation between inputs and maximal outputs (or outputs and minimal inputs), and labels any deviation from the frontier as inefficient.  For instance, the originally proposed efficiency measure in DEA is the maximum of a ratio between the weighted sum of outputs and that of inputs (see objective function of (\ref{eq:originalDEA})) and is obtained for each $\text{DMU}_k = \text{DMU}_0, k = 1, \ldots, n$, in the sample by solving the LP equivalent of the following fractional program $n$ times.
\begin{eqnarray}
\label{eq:originalDEA}
\max_{u, v}\quad\displaystyle\frac{\sum_{r = 1}^{s}{y_{r,0}u_r }}{\sum_{i = 1}^{m}{x_{i,0}v_i}} & &  \\
s.t. \quad\displaystyle\frac{\sum_{r = 1}^{s}{y_{r,k}u_r}}{\sum_{i = 1}^{m}{x_{i,k}v_i}} & \leq 1 & \quad (k = 1, \ldots, n)  \nonumber\\
u_r & \geq & 0 \quad (r = 1, \ldots, s)  \nonumber\\
v_i & \geq & 0 \quad (i = 1, \ldots, m)\nonumber
\end{eqnarray}
where $X = x_{i,k} \in \mathbb{R}^{m \times n}$ are the input parameters, $Y = y_{r, k} \in \mathbb{R}^{s \times n}$ are the output parameters, $u$ and $v$ are the variables for output and input weights respectively  \citep{CCR}.  Here, the inequality constraint is imposed to ensure that the estimated efficient frontier envelops all the sample data points.  DEA essentially generalizes the so-called productivity ratio of a single output to a single input to the case of multiple outputs and multiple inputs.

The LP equivalent of (\ref{eq:originalDEA}) is commonly known as the (primal) CCR model and is one of the three representative basic DEA models together with the BCC \citep{BCC} and additive \citep{Additive} models.  The output-oriented formulations\footnote{DEA models can be classified as  output-oriented, input-oriented or base-oriented.  While output-oriented DEA models focus on output augmentation to achieve efficiency (i.e. outputs are controllable), input-oriented DEA models are concerned with minimizing the amount of inputs required to produce a certain amount of outputs (i.e. inputs are controllable).  Base-oriented DEA models are concerned with determining the optimal mix of inputs and outputs (i.e. both inputs and outputs are controllable).  For input-oriented formulations of the CCR and BCC models presented in this paper, refer to \citep{CCR} and \citep{BCC} respectively.}  of the primal and dual CCR models are given by (\ref{eq:primalCCR}) and (\ref{eq:dualCCR}) respectively.
\begin{eqnarray}
\label{eq:primalCCR}
\max_{u, v}\quad\sum_{r = 1}^{s}{y_{r,0}u_r} & &\\
s.t.\quad\sum_{i = 1}^{m}{x_{i, 0}v_i} & = & 1\nonumber\\
\sum_{r = 1}^{s}{y_{r,k}u_r} & \leq & \sum_{i = 1}^{m}{x_{i,k}v_i} \quad (k = 1, \ldots, n)\nonumber\\
u_r & \geq & 0 \quad (r = 1, \ldots, s)\nonumber\\
v_i & \geq & 0 \quad (i = 1, \ldots, m);\nonumber
\end{eqnarray}
\begin{eqnarray}
\label{eq:dualCCR}
\max_{\theta_0, \lambda}\quad\theta_0 & &\\
s.t.\quad x_{i, 0} & \geq & \sum_{k = 1}^{n}{x_{i, k}\lambda_k} \quad (i = 1, \ldots, m)\nonumber\\\
\sum_{k = 1}^{n}{y_{j, k}\lambda_k } & \geq & y_{j, 0}\theta_0 \quad (j = 1, \ldots, s)\nonumber\\
\lambda_k & \geq & 0 \quad (k = 1, \ldots, n).\nonumber
\end{eqnarray}
The output-oriented BCC model is obtained when the above CCR model (\ref{eq:dualCCR}) is augmented by adding a convexity constraint, $\sum_{k = 1}^n \lambda_k = 1$.  This convexity constraint accounts for variable returns to scale (VRS) production technology; i.e., without this constraint, the model assumes constant returns to scale (CRS) production technology.  Basically, the CCR and BCC models differ only in their assumptions of the underlying production technology.

Both CCR and BCC models have been criticized for taking a radial approach to measuring efficiency.  In the radial approach, inputs and outputs are assumed to change proportionally, and such an assumption makes the CCR and BCC models prone to neglect non-radial input excesses and output shortfalls.  In contrast, the additive model takes a non-radial approach to efficiency estimation by directly dealing with input and output slacks and is also free of input-output orientation.  The formal definition of the dual additive model with VRS production technology\footnote{The corresponding formulation with CRS technology is a special instance without the convexity constraint $\sum_{k = 1}^n \lambda_k = 1$ in (\ref{eq:dualAdd}) (or without the variable $w$ in (\ref{eq:primalAdd})).} is given by
\begin{eqnarray}
\label{eq:dualAdd}
\max_{s^-, s^+, \lambda}\quad\ Z_0  & = & \sum_{i = 1}^{m} {s_{i, 0}^{-}}+ \sum_{r = 1}^{s} {s_{r, 0}^{+}}  \\
s.t. \quad \sum_{k = 1}^{n} \lambda_k y_{r, k}  & = & y_{r, 0} + s_{r, 0}^{+} \quad (r = 1, \ldots,s)\nonumber\\
\sum_{k = 1}^{n} \lambda_k x_{i, k}  & = & x_{i, 0} + s_{i, 0}^{-} \quad (i = 1, \ldots,m)\nonumber\\
\sum_{k =1}^{n} \lambda_{k} & = &1 \nonumber\\
\lambda_{j} & \geq & 0 \quad (j = 1, \ldots, n)\nonumber\\
s_{i, 0}^{-} & \geq & 0 \quad (i = 1, \ldots, m)\nonumber\\
s_{r, 0}^{+} & \geq & 0 \quad (r= 1, \ldots, s)\nonumber
\end{eqnarray}
where $s_{i, 0}^{-}, i = 1, \ldots, m$ and $s_{r, 0}^{+} , r = 1, \ldots, s$ represent the respective input and output slacks.  Its associated primal model is given by
\begin{eqnarray}
\label{eq:primalAdd}
\min_{u, v, w} \quad\sum_{i = 1}^{m}{x_{r,0}v_i} -\sum_{r = 1}^{s}{y_{r,0}u_r} + w & &\\
s.t. \quad\sum_{i=1}^{m}x_{i,k}v_{i} - \sum_{r=1}^{s}y_{r,k}u_{r} + w & \geq & 0, \;  k = 1,\ldots,n\nonumber\\
u_r & \geq & 1 \quad (r = 1, \ldots, s)\nonumber\\
v_i &\geq & 1 \quad (i = 1, \ldots, m).\nonumber
\end{eqnarray}
One can solve (\ref{eq:primalAdd}) for all $n$ DMUs simultaneously by solving the following LP,
\begin{eqnarray}
\label{eq:primalAddAll}
  \min_{u,v,w} && \sumTo{k}{1}{n}\left(\sumTo{i}{1}{m}x_{i,k}v_{i,k}-\sumTo{r}{1}{s}y_{r,k}u_{r,k}+w_k\right) \\
  s.t. && \sumTo{i}{1}{m}x_{i,j}v_{i,k} - \sumTo{r}{1}{s}y_{r,j}u_{r,k} + w_k \geq 0, \; j = 1,\ldots,n; k = 1,\ldots,n \nonumber\\
   && u_{r,k} \geq 1, \; r = 1,\ldots,s; k = 1,\ldots,n \nonumber\\
   && v_{i,k} \geq 1, \; i = 1,\ldots,m; k = 1,\ldots,n.\nonumber
\end{eqnarray}

In the next section, we develop a joint variable selection method for this additive model (\ref{eq:primalAddAll}).  It should, however, be noted that the proposed variable selection method can be readily adapted for various formulations of DEA models.  We should also note that unless otherwise mentioned, the same notations for variables and parameters introduced in this section are used throughout the paper.

\section{Joint Variable Selection}
\label{sec:GL}
In DEA, we are often in a situation where we want to select a small number of most relevant input variables across the DMUs.  One popular variable selection approach through convex optimization is the Lasso \citep{tib}.  It is a simple regularization technique, which adds an $l_1$-norm (sum of the absolute values) of the variables to the original objective function.  Due to the special geometric properties of the $l_1$-ball, the solution to the Lasso problem is sparse, i.e. only a small number of entries are non-zero, and these correspond to the selected variables.  Although the Lasso was originally designed for variable selection in linear regression models, it can be easily  extended to DEA models.  For instance, for the additive model (\ref{eq:primalAddAll}), its respective Lasso formulation is given by the following LP
\begin{eqnarray}
  \min_{u,v,w} && \sumTo{k}{1}{n}\left(\sumTo{i}{1}{m}x_{i,k}v_{i,k}-\sumTo{r}{1}{s}y_{r,k}u_{r,k}+w_k + \lambda\sumTo{i}{1}{m}v_{i,k}\right) \\
  s.t. && \sumTo{i}{1}{m}x_{i,j}v_{i,k} - \sumTo{r}{1}{s}y_{r,j}u_{r,k} + w_k \geq 0, \; j = 1,\ldots,n; k = 1,\ldots,n \nonumber\\
   && u_{r,k} \geq 1, \; r = 1,\ldots,s; k = 1,\ldots,n\nonumber \\
   && v_{i,k} \geq 1, \; i = 1,\ldots,m; k = 1,\ldots,n, \nonumber
\end{eqnarray}
for some $\lambda > 0$, the regularization parameter that controls the level of sparsity in the solution.  We do not need absolute values of $v$ in the objective function since they are constrained to be positive.  Note that in this case, the solution is, in fact, sparse only after a ``shift," by subtracting 1 from each entry.  We should also note that one can readily incorporate output variables into variable selection by adding $\lambda\sumTo{r}{1}{s}u_{r,k}$ to the objective function.  For an application of the Lasso on DEA models, readers can refer to \citet{SongKachaniKani}.

The results obtained from the Lasso formulation may be hard to interpret because the selection of the variables is not guaranteed to be consistent across all the DMUs.  For instance, for a given variable $i$, $v_{i,k}$ may not be selected for all $k$'s (i.e. across all DMUs).  In fact, it has been shown that the Lasso tends to select only one variable from a group of highly correlated variables and does not care which one is selected \citep{zou2005regularization}.  Hence, an approach that enforces selection consistency across the DMUs is called for.  Note that if we stack the column vectors $v_k$'s as a matrix $V$, then our goal is to select a small number of rows from $V$.

\subsection{Group Sparsity-inducing Regularization}
Before we discuss our approach to joint variable selection for DEA, we need to introduce a more general regularization technique, group Lasso \citep{yuanlin}, which is an extension to Lasso and tends to induce variable sparsity at group level, i.e. to select a small number of groups of correlated variables.  It achieves this goal via $l_{2,1}$-regularization,
\begin{equation}\label{eq:groupLasso}
    \min_\beta\;\;\ F(\beta) + \lambda\sumTo{j}{1}{J}\|\beta_j\|_2,
\end{equation}
which is the sum of the group $l_2$-norms with a pre-defined grouping of the variables $\setTo{\beta_j}{j}{1}{J}$.  The original group Lasso problem considered in \citet{yuanlin} has $F(\beta) := \frac{1}{2}\|X\beta-y\|_2^2$, i.e. least-squares regression.  The same regularization technique has also been applied to logistic regression \citep{mgb}.

The requirement of a pre-defined grouping of the variables is often a limiting factor for applications of group Lasso. However, in the case of DEA joint variable selection, the grouping structure is readily available -- we simply group the variables in $V$ by rows or DMUs.  We can then solve a special constrained version of the group Lasso with the loss function,
\begin{equation}
F(u,v,w) := \sumTo{k}{1}{n}\left(\sumTo{i}{1}{m}x_{i,k}v_{i,k}-\sumTo{r}{1}{s}y_{r,k}u_{r,k}+w_k\right).
\end{equation}
The nonzero entries in the group-sparse solution that we obtain (after shifting) are then guaranteed to be consistent across all the DMUs.

The standard group Lasso problem has been studied extensively in the machine learning and optimization literature, and a number of convex optimization algorithms (e.g. \citet*{slep, mgb, qin2010efficient, friedlander_spg, sparsa, yuanlin}) have been proposed to solve it.  However, most of these algorithms are designed to solve the unconstrained group Lasso problem.  In the subsequent sections, we propose a new tailored algorithm based on the alternating direction method of multipliers (ADMM) to solve our special constrained group Lasso problem.

\subsection{Problem Formulation}
Using a change of variables $\tilde{u} = u - e$ and $\tilde{v} = v - e$, we can transform the original additive model (\ref{eq:primalAddAll}) into the following model with zero lower bounds on $\tilde{u}$ and $\tilde{v}$:
\begin{eqnarray}
  \min_{\tilde{u},\tilde{v},w} && \sumTo{k}{1}{n}\left(\sumTo{i}{1}{m}x_{i,k}\tilde{v}_{i,k}-\sumTo{r}{1}{s}y_{r,k}\tilde{u}_{r,k} + w_k\right) \\
  s.t. && \sumTo{i}{1}{m}x_{i,j}(\tilde{v}_{i,k}+1) - \sumTo{r}{1}{s}y_{r,j}(\tilde{u}_{r,k}+1) + w_k \geq 0, \; j = 1,\cdots,n; k = 1,\cdots,n \nonumber\\
   && \tilde{u}_{r,k} \geq 0, \; r = 1,\ldots,s; k = 1,\ldots,n \nonumber\\
   && \tilde{v}_{i,k} \geq 0, \; i = 1,\ldots,m; k = 1,\ldots,n. \nonumber
\end{eqnarray}
For joint variable selection on $v$, we propose to apply the group Lasso regularization on $\tilde{v}$; i.e.
\begin{eqnarray}
  \min_{\tilde{u},\tilde{v},w} && \sumTo{k}{1}{n}\left(\sumTo{i}{1}{m}x_{i,k}\tilde{v}_{i,k}-\sumTo{r}{1}{s}y_{r,k}\tilde{u}_{r,k} + w_k\right) + \lambda\sumTo{i}{1}{m}\|\tilde{v}^i\|_2 \\
  s.t. && \sumTo{i}{1}{m}x_{i,j}(\tilde{v}_{i,k}+1) - \sumTo{r}{1}{s}y_{r,j}(\tilde{u}_{r,k}+1) + w_k \geq 0, \; j = 1,\ldots,n; k = 1,\ldots,n \nonumber\\
   && \tilde{u}_{r,k} \geq 0, \; r = 1,\ldots,s; k = 1,\ldots,n \nonumber\\
   && \tilde{v}_{i,k} \geq 0, \; i = 1,\ldots,m; k = 1,\ldots,n, \nonumber
\end{eqnarray}
where $\tilde{v}^i$ is the vector of $\tilde{v}_{i,k}$ for $k=1,\ldots,n$.
This formulation can be interpreted as selecting the elements of $v$ which are sufficiently large, i.e. larger than the lower bounds.  If joint variable selection on $u$ is also needed, we can apply the group Lasso regularization on $u$ in a similar way.  For notational simplicity, we drop the $\sim$ signs for $\tilde{u},\tilde{v},\tilde{X}$, and $\tilde{Y}$ from now on.

Writing in matrix form, the group lasso formulation for the additive model of DEA is
\begin{eqnarray}\label{eq:dea-gl}
  \min_{U,V,W} && tr(X\transp V - Y\transp U + W) + \lambda\sumTo{i}{1}{m}\|v^i\|_2 \\
  \nonumber s.t. && X\transp V - Y\transp U + W + B \geq 0 \\
  \nonumber  && U \geq 0 \\
  \nonumber  && V \geq 0,
\end{eqnarray}
where $v^i$ is $i$-th row of $V$, and $W = \big( ew_1 \; \cdots \; ew_n \big)$.  Introducing non-negative slack variables, we can transform \eqref{eq:dea-gl} into
\begin{eqnarray}\label{eq:dea-gl-slack}
  \min_{U,V,\barV,W,S} && tr(X\transp V - Y\transp U + W) + \lambda R(\barV) \\
  \nonumber s.t. && X\transp V - Y\transp U + W + B = S_x, \; S_x \geq 0 \\
  \nonumber  && U = S_u, \; S_u \geq 0 \\
  \nonumber  && V = \barV = S_v, \; S_v \geq 0,
\end{eqnarray}
where $R(\barV) = \sumTo{i}{1}{m}\|v^i\|_2$.
For the ease of visualization, we write the problem in terms of vectorized decision variables, i.e. stacking columns on top of each other, which is the same as the (:) operator in Matlab.  We use the lower-case letter to denote the vectorized version of the matrix counterpart, which is denoted by upper-case.  The following are some elementary transformations:
\begin{equation}
    X\transp V \; \rightarrow \; \left(
                                   \begin{array}{ccc}
                                     X\transp &  &  \\
                                      & \ddots &  \\
                                      &  & X\transp \\
                                   \end{array}
                                 \right) v = \barX v \; , \; tr(X\transp V) \;\rightarrow\; x\transp v \; , \; W \;\rightarrow\; \left(
                                                                                                                                   \begin{array}{ccc}
                                                                                                                                     e &  &  \\
                                                                                                                                      & \ddots &  \\
                                                                                                                                      &  & e \\
                                                                                                                                   \end{array}
                                                                                                                                 \right) w = \barI w.
\end{equation}
Problem \eqref{eq:dea-gl-slack} then becomes
\begin{eqnarray}
  \min_{s\geq 0, v,u,w,\barv} && x\transp v - y\transp u + e\transp w + \lambda R(\barv) \\
  \nonumber s.t. && \barX v - \barY u + \barI w + b = s_x \\
  \nonumber && v = s_v \\
  \nonumber  && u = s_u \\
  \nonumber && v = \barv.
\end{eqnarray}
By writing the above problem in a more compact form, we obtain
\begin{eqnarray}\label{eq:dea-gl-constr}
  \min_{s\geq 0, z,\barv} && c\transp z + \lambda R(\barv) \\
  \nonumber s.t. && \left(
                      \begin{array}{c}
                        A_s \\
                        A_v \\
                      \end{array}
                    \right) z  + \left(
                      \begin{array}{c}
                        b \\
                        0 \\
                      \end{array}
                    \right)= \left(
                                  \begin{array}{c}
                                    s \\
                                    \barv \\
                                  \end{array}
                                \right),
\end{eqnarray}
where $z = \left(
             \begin{array}{c}
               v \\
               u \\
               w \\
             \end{array}
           \right), \; c = \left(
                             \begin{array}{c}
                               x \\
                               -y \\
                               e \\
                             \end{array}
                           \right), \; A_s = \left(
                                               \begin{array}{ccc}
                                                 \barX & -\barY & \barI \\
                                                 I & 0 & 0 \\
                                                 0 & I & 0 \\
                                               \end{array}
                                             \right), \; A_v = \left(
                                                                 \begin{array}{ccc}
                                                                   I & 0 & 0 \\
                                                                 \end{array}
                                                               \right), \; s = \left(
                                                                                 \begin{array}{c}
                                                                                   s_x \\
                                                                                   s_v \\
                                                                                   s_u \\
                                                                                 \end{array}
                                                                               \right)
                                                               .
$

We note that the BCC model (for all $n$ DMUs) with the group Lasso regularization is
\begin{eqnarray}
\label{eq:bccLasso}
  \min_{u,v,w} && \sumTo{k}{1}{n}\left(x_k\transp v_k + w_k\right) + \lambda\sumTo{i}{1}{m}\|v^i\|_2 \\
  s.t. && y_k\transp u_k = 1, \; k = 1,\ldots,n \nonumber\\
   && X\transp v_k - Y\transp u_k + w_k e\transp + b \geq 0, \; k = 1,\ldots,n \nonumber\\
   && u_k \geq 0, \; k = 1,\ldots,n \nonumber\\
   && v_k \geq 0, \; k = 1,\ldots,n,\nonumber
\end{eqnarray}
where $v_k$ is the vector of $v_{i,k}, i=1,\ldots,m$, and similar definitions apply to $u_k, x_k$, and $y_k$.  $X$ and $Y$ are matrices whose columns are $x_k$'s and $y_k$'s respectively.  The corresponding problem for the CCR model is the same except that there is no variable $w$.  Similar transformation involving non-negative slacks applied to the additive model can also be applied to the above model (\ref{eq:bccLasso}).  Likewise, the optimization algorithm described in the next section can be used to solve this problem as well.

\subsection{Optimization Algorithm}
The alternating direction method of multipliers (ADMM) was first proposed in the 1970s \citep{gabay1976dual, glowinski1975adal}.  It belongs to the family of the classical augmented Lagrangian (AL) method \citep{powell1972nonlinear, rockafellar1973multiplier, hestenes1969multiplier}, which iteratively solves the linearly constrained problem
\begin{eqnarray}\label{eq:lin_constr_prob}
  \min_x && F(x) \\
  \nonumber s.t. && Ax = b.
\end{eqnarray}
The augmented Lagrangian of problem \eqref{eq:lin_constr_prob} is $\mathcal{L}(x,\gamma) = F(x) + \gamma^T(b-Ax) + \frac{1}{2\mu}\|Ax-b\|^2$, where $\gamma$ is the Lagrange multiplier and $\mu$ is the penalty parameter for the quadratic infeasibility term.  The AL method minimizes $\mathcal{L}(x,\gamma)$ followed by an update to $\gamma$ in each iteration.

For a structured unconstrained problem
\begin{equation}\label{eq:struct_uncon_prob}
    \min_x F(x) \equiv f(x) + g(Ax),
\end{equation}
where both functions $f(\cdot)$ and $g(\cdot)$ are convex, we can decouple the two functions by introducing an auxiliary variable $y$ and transform problem \eqref{eq:struct_uncon_prob} into an equivalent linearly constrained problem
\begin{eqnarray}\label{eq:struct_lin_con}
  \min_{x,y} && f(x) + g(y) \\
  \nonumber s.t. && Ax = y.
\end{eqnarray}
The augmented Lagrangian of this problem is
\begin{equation}
    \mathcal{L}(x,y,\gamma) = f(x) + g(y) + \gamma^T(y-Ax) + \frac{1}{2\mu}\|Ax-y\|^2.
\end{equation}
ADMM (Algorithm \ref{alg:adal}) finds the approximate minimizer of $\mathcal{L}(x,y,\gamma)$ by alternatively optimizing with respect to $x$ and $y$ once.  This is often desirable because joint minimization of $\mathcal{L}(x,y,\gamma)$ even approximately could be hard.
\begin{algorithm}
\caption{ADMM}
\begin{algorithmic}[1]\label{alg:adal}
\STATE Choose $\xsupk{\gamma}{0}$.
\FOR{$k = 0,1,\ldots,K$}
    \STATE $\xkpone \gets \arg\min_{x}\mathcal{L}(x, \xsupk{y}{k}, \gammak )$
    \STATE $\ykpone \gets \arg\min_{y}\mathcal{L}(\xkpone, y, \gammak )$
    \STATE $\gammakpone \gets \gammak - \frac{1}{\mu}(A\xkpone - \ykpone)$
\ENDFOR
\RETURN $\xsupk{y}{K}$
\end{algorithmic}
\end{algorithm}

Our strategy is to apply Algorithm \ref{alg:adal} to solve problem \eqref{eq:dea-gl-constr}.  First, we write down the augmented Lagrangian of the problem,
\begin{equation}\label{eq:dea-augLag}
    \calL(z,s,\barv,\gamma_s,\gamma_v) := c\transp z + \lambda R(\barv) - \gamma_s\transp(A_s z + b -s) + \frac{1}{2\mu}\|A_s z + b -s\|^2 - \gamma_v\transp(A_v z  -\barv) + \frac{1}{2\mu}\|A_v z  -\barv\|^2.
\end{equation}
Next, we minimize with respect to $z,s,\barv$ sequentially.  The subproblem with respect to $z$ can be simplified to solving a sparse linear system with a fixed left-hand-side,
\begin{equation}\label{eq:z-subprob}
    \frac{1}{\mu}(A_s\transp A_s + A_v\transp A_v)z = A_s\transp\gamma_s + A_v\transp\gamma_v - c + \frac{1}{\mu}(A_s\transp s + A_v\transp \barv).
\end{equation}
As long as we keep $\mu$ constant, we can compute the Cholesky factor of the left-hand-side for once and cache it for subsequent iterations, where the computation for this step is almost as cheap as a gradient step (via forward/backward substitution).
The subproblem with respect to $s$ is a projection problem onto the non-negative orthant,
\begin{equation}\label{eq:s-subprob}
    \min_{s\geq 0} \; \frac{1}{2}\|A_s z + b - s - \mu\gamma_s\|^2.
\end{equation}
We can obtain the solution easily by $(A_s z + b - \mu\gamma_s)^+$, where $(\cdot)^+$ is an element-wise truncation operation at 0.
The subproblem with respect to $\barv$ is the proximal problem associated with the group lasso penalty,
\begin{equation}\label{eq:vbar-subprob}
    \min_{\barv} \; \frac{1}{2}\|A_v z  - \barv - \mu\gamma_v\|^2 + \mu\lambda R(\barv).
\end{equation}
The optimal solution can be computed in closed-form: $\barv^* = \calT_{\mu\lambda}(A_v z  - \mu\gamma_v)$, where $\calT$ is the block soft-thresholding operator such that the $i$-th block of $\barv^*$, $[\barv^*]_i = \frac{A_v z  - \mu\gamma_v}{\|A_v z  - \mu\gamma_v\|_2}\max(0,\|A_v z  - \mu\gamma_v\|_2 - \mu\lambda)$, for $i=1,\ldots,m$.

\subsection{Convergence}
The convergence of ADMM has been established for the case of two-way splitting as above.  We restate the results from \citep{eckstein1992douglas} in the following theorem.
\begin{thm}\label{thm:admm_conv}
Consider problem \eqref{eq:struct_lin_con}, where both $f$ and $g$ are proper, closed, convex functions, and $A \in \mathbb{R}^{n\times m}$ has full column rank.  Then, starting with an arbitrary $\mu > 0$ and $x^0,y^0 \in \mathbb{R}^m$, the sequence $\{x^k,y^k,\gamma^k\}$ generated by Algorithm \ref{alg:adal} converges to a Kuhn-Tucker pair $\big((x^*,y^*),\gamma^*\big)$ of problem \eqref{eq:struct_lin_con}, if \eqref{eq:struct_lin_con} has one.  If \eqref{eq:struct_lin_con} does not have an optimal solution, then at least one of the sequences $\{(x^k,y^k)\}$ and $\{\gamma^k\}$ diverges.
\end{thm}
It is known that $\mu$ does not have to decrease to a very small value (or can simply stay constant) in order for the method to converge to the optimal solution of problem \eqref{eq:struct_lin_con} \citep{nocedal, bertsekas1999nonlinear}.

We observe that problem \eqref{eq:dea-gl-constr} can be treated as a two-way splitting, with variables $z$ and $\left(
                                  \begin{array}{c}
                                    s \\
                                    \barv \\
                                  \end{array}
                                \right)$, and obviously the matrix $\left(
                      \begin{array}{c}
                        A_s \\
                        A_v \\
                      \end{array}
                    \right)$ has full column rank.  Hence, Theorem \ref{thm:admm_conv} applies to our ADMM algorithm.

In the next section, we introduce benchmark variable selection methods, against which we evaluate the performance of the proposed joint variable selection method.  We refer to our group Lasso-based variable selection method as GL method now on.

\section{Benchmarks}
\label{sec:benchmark}
Among prevailing variable selection methods in the DEA literature,\footnote{Readers can refer to \citet{JohnsonNataraja} for detailed reviews of prevailing variable selection methods in the DEA literature.} the four most widely-used approaches are 1) the regression-based (RB) test \citep{Ruggiero}, 2) the efficiency contribution measure (ECM) method \citep{PastorEtAl}, 3) the principal component analysis (PCA-DEA) \citep{UedaHoshiai, AdlerGolany}, and 4) the bootstrapping method \citep{SimarWilson}.  Nataraja and Johnson \citeyearpar{JohnsonNataraja} have already evaluated these four approaches and reported their performance in their study.  According to their results, RB test and ECM method are best suited for various sample sizes provided that there is low correlation among variables while PCA-DEA and bootstrapping method show some limitations.  The two major limitations of PCA-DEA are: first, as it replaces the original variables with principal components (PCs), the original data set is not retained, and therefore it is impossible to recover true efficiency levels, and second, it is vulnerable to the curse of dimensionality.  The vital issue with the bootstrapping method is that it involves the heavy computational burden, and yet, has the weakest performance among the four.  Consequently, we select ECM method and RB test as the benchmark variable selection methods, against which the performance of our proposed method is measured.

\subsection{Efficiency Contribution Measure (ECM) Method}
The efficiency contribution measure (ECM) method evaluates the effect of a candidate variable $x_{\text{cand}}$ on efficiency computation by comparing two DEA formulations -- one with the candidate variable and one without it.  The ECM of $x_{\text{cand}}$ for a particular $\text{DMU}_0$, denoted by $\gamma_0$, is a single scalar measure that quantifies the marginal impact of $x_{\text{cand}}$ on the measurement of efficiency.  In essence, the ECM method performs a statistical test to determine the statistical significance of $x_{\text{cand}}$'s contribution when measured by means of ECMs.  It should be noted that the ECM method consists of two procedures for the progressive selection of variables -- a forward selection (addition of variables) and backward elimination (removal of variables) and only supports radial DEA models, such as the CCR and BCC models.

To provide further technical details, suppose $\gamma = (\gamma_1, ..., \gamma_n)$ are the observed ECMs of a random sample, $\Gamma = (\Gamma_1, ..., \Gamma_n)$ drawn from a population $(\Gamma, F)$  where $\Gamma$ being a random variable distributed according to $F$, a cumulative density function on $[1, \infty)$.  The underlying idea of the ECM method is that if $x_{\text{cand}}$ is an irrelevant variable, then the impact it has on efficiency evaluation should be negligible, and high values of $\Gamma$ associated with $x_{\text{cand}}$ are unlikely to be observed.  For the statistical test, two additional parameters $\bar{\gamma}$ and $p_0$ are introduced.  $\bar{\gamma}$ represents the tolerance level for the degree of efficiency score change caused by $x_{\text{cand}}$, and $p_0$ represents the tolerance level for the proportion of DMUs whose associated efficiency score change exceeds $\bar{\gamma}$.   $x_{\text{cand}}$ is considered relevant to the production process if more than $p_0$\% of DMUs have associated efficiency score change greater than $\bar{\gamma}$.  More formally, a hypothesis test with a binomial test statistic is performed to see if the marginal impact of this candidate variable on efficiency estimation is significant.  For technical details and applications of the ECM method, readers can refer to \citet{PastorEtAl} and \citet{ChenJohnson} respectively.

\subsection{Regression-based (RB) Test}
In the regression-based (RB) test, initial efficiency estimates obtained from the set of known production variables are regressed against the set of candidate variables.  The formal representation of the regression model is given by
\begin{equation}
E = \alpha + \beta_2 x_2 + \beta_3 x_3 + \cdots + \beta_m x_m + \epsilon,
\end{equation}
where $E$ is the efficiency score obtained from the DEA model including only an output variable $y$ and without loss of generality, an input variable $x_1$, and $x_2$ through $x_m$ are candidate variables.  If the coefficient $\beta_i$ in the regression is statistically significant at a given level of significance and has the proper sign (i.e. $\beta > 0$ for input variables and $\beta < 0$ for output variables), the candidate variable $x_i$ is considered relevant to the production process and is added to the DEA model.  A new efficiency score $E$ is then computed with this updated DEA model, and the new candidate set is tested.  This process is repeated until all candidate variables are either found irrelevant or included in the model, and there are no more remaining variables to be tested.  For technical details of the RB test, readers can refer to Ruggiero \citeyearpar{Ruggiero}.

\section{Experimental Design and Data Generation}
\label{sec:mc}
In our simulation study, we focus on output-oriented radial DEA models, namely the CCR \citep{CCR} and BCC \citep{BCC} models.  The reason behind this particular choice of models is that one of the benchmark methods, the ECM method, is not compatible with non-radial models.  A presentation based on input-oriented formulations or non-radial models may be similarly developed.

In practical applications of DEA, the true form of production process is mostly unknown, and the observed data used for estimating this unknown production function are often limited and contain measurement errors.  These are the major setbacks in evaluating the practical importance of theoretical results in DEA, as a matter of fact.  In order to overcome this problem, this study uses Monte Carlo simulations to generate a large number of observations for a plausible production process, the form of which is known.

The production process we consider is the linearly homogeneous Cobb-Douglas function, in which a number of inputs, represented by vector $x$, are used to produce a single output $y$; i.e.
\begin{equation}
\label{eq:cobb_douglas}
y_{1, k} = \beta\prod_{i = 1}^m x_{i, k}^{\alpha^i}, \quad k = 1,\ldots, n
\end{equation}
where $\alpha_i$ and $\beta$ are assumed to be known parameters.  The parameter  $\alpha_i$ here plays an important role in the production model.  First, it defines the returns to scale (RTS) specification for the production process; i.e. $\sum_{i = 1}^m \alpha_i = 1$ indicates a CRS production process while $\sum_{i = 1}^m \alpha_i < 1$ indicates a VRS production process.  Second, mathematically, $\alpha_i$ indicates the flexibility of production with respect to the input $x_i$ under efficient production.  Intuitively speaking, it represents the importance of the input $x_i$ in the production process other things being equal.  For our simulation study, we set $\alpha_i = 1/m$ where $m$ is a predetermined value as the base case scenario and $\beta = 1$ throughout.

By incorporating a technical inefficiency of the $k^{th}$ DMU denoted by $\epsilon_k \in [0, 1]$ into (\ref{eq:cobb_douglas}), we obtain
\begin{equation}
\label{eq:cobb_douglas2}
y_{1, k} = \beta\prod_{i = 1}^m x_{i, k}^{\alpha^i}\epsilon_k, \quad k = 1,\ldots, n.
\end{equation}
An alternative additive representation of (\ref{eq:cobb_douglas2}) is
 \begin{equation}
\label{eq:cobb_douglas3}
Y_{1, k} = B + \sum_{i = 1}^m\alpha_i X_{i,k} - u_k, \quad k = 1,\ldots, n
\end{equation}
where $Y_{1, k} = \ln y_{1, k}, B = \ln\beta, X_{i, k} = \ln x_{i, k}$, and $u_k = \ln\epsilon_k$.  The efficiency component $u_k \geq 0$ represents the shortfall of output from the production frontier.

With respect to statistical distributions of variables $u$ and $X$, consistent with previous studies \citep{Smith, JohnsonNataraja}, $u$ is drawn from a half-normal distribution with mean zero and variance $\sigma^2$, which we vary to obtain a sample average efficiency score of 85\% and is assumed to be uncorrelated with any $X$ in order to generate a realistic range of inefficiency values.  The values of $X$ are generated from a uniform distribution on an interval $[a, b]$ and are exponentiated and used in $(\ref{eq:cobb_douglas2})$ to yield the values of $y$ and $x$ to be used in DEA models.

Smith's \citeyearpar{Smith} simulation study involving a Cobb-Douglas production function has shown that the performance of DEA in estimating true efficiencies diminishes as the number of inputs in the production process increases.  We therefore assume that the true production process is determined by three inputs, $x_1$, $x_2$, and $x_3$ only.  In the base case scenario, we also independently generate an irrelevant random variable $x_4$ from a uniform distribution on the same interval $[a, b]$ to maintain symmetry with the three other relevant inputs.
\begin{table}[htbp]
\begin{center}
\caption{\footnotesize Outline of the Experimental Scenarios}
\begin{minipage}[t]{36pc}
\scriptsize Table~\ref{tb_testcases} delineates the experiments used for evaluating the performance of the three variable selection methods -- GL method, ECM method and RB test.  A total of 12 experimental scenarios are considered for each CRS and VRS production process.  The respective values of input contribution parameter $\alpha$ for CRS and VRS production processes are shown on the third column separated by a semicolon.
\begin{center}
{\footnotesize\tabcolsep=5pt\begin{tabular*}{36pc}{@{\extracolsep{\fill}} c l l l}
\hline%
\multicolumn{1}{l}{Experiment} & \multicolumn{1}{l} {Correlation} & \multicolumn{1}{l} {Input Contribution} & \multicolumn{1}{l} {Description} \\
\multicolumn{1}{l}{} & \multicolumn{1}{l} {Between Inputs} & \multicolumn{1}{l} {to Output} & \multicolumn{1}{l} {} \\
\hline%
1 & Independently generated & $\alpha_i = 1/3$; $\alpha_i = 1/4, i = 1, 2, 3$  & Base case\\
2 & $\rho_{1, 2} = 0.8, \rho_{1, 3} = 0.2$ & $\alpha_i = 1/3$; $\alpha_i = 1/4, i = 1, 2, 3$  & Correlated inputs\\
3 & $\rho_{1, 2} = \rho_{1, 3} = 0.8$& $\alpha_i = 1/3$; $\alpha_i = 1/4, i = 1, 2, 3$  & Highly correlated inputs\\
4 & Independently generated & $\alpha_1 = 1/3, \alpha_2 = 4/9, \alpha_3 = 2/9$;  & Input contribution to \\
& & $\alpha_1 = 1/4, \alpha_2 = 1/3, \alpha_3 = 1/6$&output varied\\
5 & $\rho_{1, 2} = 0.8, \rho_{1, 3} = 0.2$& $\alpha_1 = 1/3, \alpha_2 = 4/9, \alpha_3 = 2/9$;  & Correlated inputs and input \\
& & $\alpha_1 = 1/4, \alpha_2 = 1/3, \alpha_3 = 1/6$&contribution to output varied\\
6 & $\rho_{1, 2} = 0.8, \rho_{1, 3} = 0.2$ & $\alpha_1 = 1/3, \alpha_2 = 2/9, \alpha_3 = 4/9$;  & Correlated inputs and input\\
& & $\alpha_1 = 1/4, \alpha_2 = 1/6, \alpha_3 = 1/3$& contribution to output varied\\
7 & $\rho_{1, 4} = 0.8$ & $\alpha_i = 1/3$; $\alpha_i = 1/4, i = 1, 2, 3$  & Correlated input and a \\
& & &random variable\\
8 & Independently generated& $\alpha_i = 1/3$; $\alpha_i = 1/4, i = 1, 2, 3$  & Small sample size, $n = 25$\\
9 & Independently generated& $\alpha_i = 1/3$; $\alpha_i = 1/4, i = 1, 2, 3$  & Large sample size, $n = 300$\\
10 & Independently generated& $\alpha_i = 1/4$; $\alpha_i = 1/5, i = 1, 2, 3, 4$  & Base case with one more\\
& & & relevant input $x_4$\\
11 & Independently generated& $\alpha_i = 1/2$; $\alpha_i = 1/3, i = 1, 2$  & Base case without a relevant\\
& & &  input $x_3$\\
12 & Independently generated& $\alpha_i = 1/3$; $\alpha_i = 1/4, i = 1, 2, 3$  & Base case with three irrelevant\\
& & &  inputs $x_4, x_5$ and $x_6$\\
\hline%
\end{tabular*}}%
\label{tb_testcases}
\end{center}
\end{minipage}
\end{center}
\end{table}%
We test the basic variable set consisting of $y, x_1, x_2, x_3$, and $x_4$ using the three variable selection methods -- GL method, ECM method, and RB test-- to determine the model specification.  We should note that in the remaining of the paper, ``candidate variables" refer to $x_1$, $x_2$, $x_3$, and $x_4$, ``model efficiency estimates" refer to efficiency estimates obtained from a DEA model with a set of input variables identified by each variable selection method, and ``true efficiency estimates" refer to those obtained from a DEA model with true input variables.

In addition to the base case, diversified experimental scenarios are considered.  These include varying: the covariance structure of inputs, sample size, contribution of each input to output and the dimensionality of the production process.  For those experiments concerned with correlated input variables, we adopt the following equation from Wang and Schmidt \citeyearpar{WangSchmidt} to establish the desired covariance structure of inputs.
\begin{equation}
x_i = \rho_{i, j} x_j + w\sqrt{1-\rho_{i, j}^2}, \quad i = 2, 3, 4, \quad j = 1, 2, 3, \quad i \neq j.
\end{equation}
Here, $\rho_{i, j}$ is the correlation between $x_i$ and $x_j$, and $w$ is a random variable generated from a uniform distribution on the interval $[a, b]$.  Table~\ref{tb_testcases} delineates the experimental scenarios considered in the simulation study.  We should note that these scenarios are in line with the experiments used in Nataraja and Johnson's study \citeyearpar{JohnsonNataraja}.  Each experiment is tried 100 times, and the simulation results averaged over 100 trials are presented in the next section.  It should be also noted that we tuned the parameters of the algorithms used in the GL method on the training data (10\% of the full data), and the reported results are out-of-sample results.
\begin{table}[htbp]
\begin{center}
\caption{\footnotesize Parameter Specification for the Simulation Study}
{\footnotesize\tabcolsep=0pt\begin{tabular*}{20pc}{@{\extracolsep{\fill}}l l}
\hline%
\multicolumn{1}{l}{Algorithm}& \multicolumn{1}{l}{Parameter value}\\
\hline%
ECM & $p_0 = 0.15$, $\bar{\rho} = 1.10$, $\alpha = 0.05$\\
RB & $\alpha = 0.90$\\
\hline%
\end{tabular*}}%
\label{tb_params}
\end{center}
\end{table}

Table~\ref{tb_params} presents the parameter specification for the benchmark methods, the ECM method and the RB test.  We keep these parameter values the same throughout Monte Carlo simulations.  For the ECM method, following the recommendations of Pastor et al. \citeyearpar{PastorEtAl}, we set $p_0 = 15\%$, $\bar{\gamma} = 10\%$ and the significance level $\alpha$ to $5\%$ for the hypothesis test.  Also, for comparative purposes, we use the backward procedure, which begins with the full model and then eliminates one variable that has the least impact on the efficiency calculation at each successive step.  For the RB test, we set the significance level $\alpha $ to $90\%$ following Ruggiero's \citeyearpar{Ruggiero} suggestion, and without loss of generality, we choose $x_1$ as the first variable to be included in the initial efficiency estimation assuming no prior knowledge of production input variables.

\begin{table}[htbp]
\begin{center}
\caption{\footnotesize Performance of the Variable Selection Methods for a CRS Production Process}
\begin{minipage}[t]{36pc}
\scriptsize Table~\ref{tab_perfMetric} presents MSE, Pearson's correlation coefficient and Spearman's rank correlation coefficient between the true and model efficiency estimates.  True efficiency estimates are obtained using a CCR model with true input variables while three sets of model efficiency estimates are obtained using a CCR model with GL method, ECM method and RB test selected input variables respectively.
\begin{center}
{\footnotesize\tabcolsep=5pt\begin{tabular*}{36pc}{@{\extracolsep{\fill}} c | c c c | c c c | c c c}
\hline%
\multicolumn{1}{c}{Metrics} & \multicolumn{3}{c} {MSE} &\multicolumn{3}{c}{Correlation Coefficient}&\multicolumn{3}{c}{Rank Correlation Coefficient}\\
\hline%
Experiments & GL & ECM & RB & GL & ECM & RB & GL & ECM & RB\\
\hline%
1&0.0000&0.0000&0.0118&1.0000&1.0000&0.9302&0.9998&1.0000&0.9242\\
2&0.0001&0.0013&0.0022&0.9985&0.9763&0.9789&0.9976&0.9662&0.9744\\
3&0.0001&0.0026&0.0069&0.9981&0.9599&0.9398&0.9976&0.9486&0.9288\\
4&0.0005&0.0001&0.0054&0.9908&0.9980&0.9631&0.9884&0.9975&0.9583\\
5&0.0007&0.0023&0.0027&0.9900&0.9629&0.9712&0.9857&0.9496&0.9641\\
6&0.0003&0.0008&0.0019&0.9936&0.9852&0.9874&0.9908&0.9775&0.9856\\
7&0.0003&0.0000&0.0573&0.9944&0.9988&0.6613&0.9918&0.9984&0.6356\\
8&0.0014&0.0010&0.0176&0.9660&0.9673&0.8432&0.9428&0.9445&0.8092\\
9&0.0000&0.0000&0.0100&1.0000&1.0000&0.9494&0.9999&1.0000&0.9466\\
10&0.0001&0.0005&0.0096&0.9988&0.9883&0.9369&0.9976&0.9833&0.9278\\
11&0.0000&0.0000&0.0057&1.0000&1.0000&0.9625&1.0000&1.0000&0.9583\\
12&0.0000&0.0000&0.0003&1.0000&1.0000&0.9921&0.9999&1.0000&0.9880\\
\hline
\end{tabular*}}%
\label{tab_perfMetric}
\end{center}
\end{minipage}
\end{center}
\end{table}%
\section{Results}
\label{sec:results}
Performance criteria used for evaluating the three variable selection methods -- GL method, ECM method and RB test -- can be broadly divided into three sets.  The first and second measurement criteria we consider are the mean squared error (MSE) and the correlations between the true and model efficiency estimates.  For correlation metrics, we use Pearson's and Spearman's rank correlation coefficients.  The last set consists of two measures: 1) the percentage of all DMUs correctly identified as efficient or inefficient and 2) the percentage of efficient DMUs correctly identified as efficient.  All the methods under evaluation are implemented in Matlab.  The source code for the GL method is available upon request.

The results for the first and second performance criteria for CRS and VRS production processes are presented in Table~\ref{tab_perfMetric} and \ref{tab_perfMetric2} respectively.  Table~\ref{tab_identification} provides the results for the last set of performance criteria.  We will discuss these results in relation to variations in the covariance structure of inputs, sample size, importance of inputs in the production process, and the dimensionality of the production space.
\begin{table}[htbp!]
\begin{center}
\caption{\footnotesize Performance of the Variable Selection Methods for a VRS Production Process}
\begin{minipage}[t]{36pc}
\scriptsize Table~\ref{tab_perfMetric2} presents MSE, Pearson's correlation coefficient and Spearman's rank correlation coefficient between the true and model efficiency estimates.  True efficiency estimates are obtained using a BCC model with true input variables while three sets of model efficiency estimates are obtained using a BCC model with GL method, ECM method and RB test selected input variables respectively.
\begin{center}
{\footnotesize\tabcolsep=5pt\begin{tabular*}{36pc}{@{\extracolsep{\fill}} c | c c c | c c c | c c c}
\hline%
\multicolumn{1}{c}{Metrics} & \multicolumn{3}{c} {MSE} &\multicolumn{3}{c}{Correlation Coefficient}&\multicolumn{3}{c}{Rank Correlation Coefficient}\\
\hline%
Experiment & GL & ECM & RB & GL & ECM & RB & GL & ECM & RB\\
\hline%
1&0.0002&0.0010&0.0022&0.9947&0.9768&0.9662&0.9885&0.9681&0.9565\\
2&0.0002&0.0024&0.0020&0.9943&0.9459&0.9529&0.9895&0.9242&0.9338\\
3&0.0006&0.0032&0.0028&0.9830&0.9345&0.9437&0.9757&0.9118&0.9223\\
4&0.0012&0.0030&0.0013&0.9687&0.9284&0.9687&0.9540&0.8989&0.9533\\
5&0.0009&0.0042&0.0037&0.9764&0.9153&0.9197&0.9656&0.8836&0.8879\\
6&0.0005&0.0013&0.0012&0.9862&0.9672&0.9660&0.9757&0.9471&0.9472\\
7&0.0002&0.0019&0.0056&0.9945&0.9557&0.9193&0.9833&0.9355&0.8986\\
8&0.0038&0.0046&0.0091&0.8845&0.8444&0.7842&0.7848&0.7636&0.6881\\
9&0.0001&0.0159&0.0009&0.9980&0.8562&0.9892&0.9964&0.8273&0.9860\\
10&0.0004&0.0085&0.0039&0.9871&0.8161&0.9153&0.9695&0.7492&0.8735\\
11&0.0000&0.0001&0.0008&1.0000&0.9974&0.9863&0.9994&0.9954&0.9825\\
12&0.0002&0.0014&0.0020&0.9941&0.9672&0.9527&0.9879&0.9564&0.9238\\
\hline
\end{tabular*}}%
\label{tab_perfMetric2}
\end{center}
\end{minipage}
\end{center}
\end{table}%
\subsection{Impact of Variations in the Covariance Structure of Inputs}
\begin{figure}[h!]
	\centering
	\includegraphics[height = 3.25in]{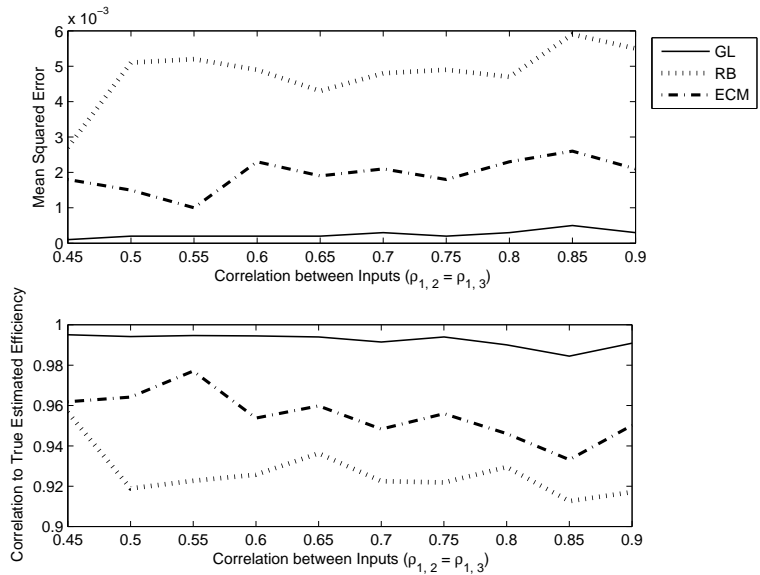}
	\caption{\footnotesize Impact of Variations in the Covariance Structure of Inputs (for a VRS Production Process)}
	\label{covstructure}
\end{figure}

In practice, input variables are often highly correlated with each other as they are all related to the scale and types of operations of DMUs being evaluated.  Accordingly, it is important to test the robustness of each variable selection method with respect to variations in the covariance structure of inputs.  Figure~\ref{covstructure} illustrates the impact of varying input correlations on the performance of the three variable selection methods for a VRS production process.  As $\rho_{1,2}$ and $\rho_{1,3}$ are varied from 0.45 to 0.90, the GL method exhibits consistently strong performance with low MSE and high correlation coefficient between the true and model efficiency estimates.   In contrast, both ECM method and RB test show fluctuations in their performance.
\begin{table}[htbp]
\begin{center}
\caption{\footnotesize Identification of Efficient and Inefficient DMUs}
\begin{minipage}[t]{36pc}
\scriptsize Table~\ref{tab_identification} presents 1) the percentage of DMUs correctly identified as efficient or inefficient and 2) the percentage of the efficient DMUs correctly identified as efficient by the three methods -- GL method, ECM method, and RB test for both CRS and VRS production frontier.
\begin{center}
{\footnotesize\tabcolsep=5pt\begin{tabular*}{36pc}{@{\extracolsep{\fill}} c | c c c | c c c | c c c | c c c}
\hline
\multicolumn{1}{c}{} & \multicolumn{6}{c}{\% Correctly Identified as Efficient/Inefficient} &  \multicolumn{6}{c}{\% of Efficient DMUs Identified as Efficient}\\
\hline%
\multicolumn{1}{c}{RTS:} & \multicolumn{3}{c}{CRS} &  \multicolumn{3}{c}{VRS}& \multicolumn{3}{c}{CRS} &  \multicolumn{3}{c}{VRS}\\
\hline%
Experiment & GL & ECM & RB & GL & ECM & RB  & GL & ECM & RB & GL & ECM & RB \\
\hline
1&100\%&100\%&98\%&99\%&97\%&97\% &100\%&100\%&82\%&98\%&89\%&89\%\\
2&100\%&95\%&99\%&99\%&92\%&94\%&97\%&46\%&88\%&96\%&57\%&67\%\\
3&100\%&95\%&94\%&98\%&91\%&91\%&96\%&30\%&28\%&84\%&39\%&41\%\\
4&99\%&100\%&98\%&96\%&91\%&96\%&92\%&98\%&85\%&87\%&65\%&86\%\\
5&99\%&95\%&98\%&97\%&88\%&89\%&86\%&41\%&74\%&84\%&36\%&47\%\\
6&98\%&95\%&99\%&97\%&94\%&95\%&81\%&50\%&95\%&86\%&67\%&75\%\\
7&99\%&100\%&89\%&99\%&94\%&93\%&98\%&98\%&9\%&99\%&80\%&72\%\\
8&96\%&96\%&88\%&90\%&85\%&78\%&89\%&88\%&57\%&82\%&75\%&58\%\\
9&100\%&100\%&99\%&100\%&92\%&99\%&100\%&100\%&87\%&100\%&46\%&95\%\\
10&100\%&98\%&96\%&99\%&81\%&90\%&99\%&91\%&79\%&99\%&47\%&74\%\\
11&100\%&100\%&99\%&100\%&100\%&99\%&100\%&100\%&89\%&100\%&100\%&95\%\\
12&100\%&100\%&98\%&99\%&97\%&95\%&100\%&100\%&100\%&100\%&87\%&89\%\\
\hline
\end{tabular*}}%
\label{tab_identification}
\end{center}
\end{minipage}
\end{center}
\end{table}%
Similarly, for a CRS production process with highly correlated input variables (i.e. Experiment 3 where $\rho_{1, 2} = \rho_{1, 3} = 0.80$), the GL method outperforms both of its benchmarks.  For instance, the GL method identifies  96\% of efficient DMUs correctly while the ECM method and RB test identify less than 30\% of them correctly.  The respective MSEs of the ECM method and RB test are also 69 and 25 times higher than that of the GL method.  When the relevant input variable $x_1$ is highly correlated with the irrelevant input variable $x_4$ in Experiment 7, the GL and ECM methods show comparably strong performance for a CRS production process.  The RB test, however, tends to choose $x_4$ as a relevant variable and has contrastingly weak performance.  For a VRS production process, the GL method outperforms both ECM method and RB test.  Overall, the GL method is most robust to variations in the covariance structure of inputs.  Consistent with previous findings, all three methods generally perform better when there is low correlation among input variables.

\subsection{Impact of Variations in Sample Size ($n$)}
\begin{figure}[h!]
	\centering
	\includegraphics[height = 2.5in]{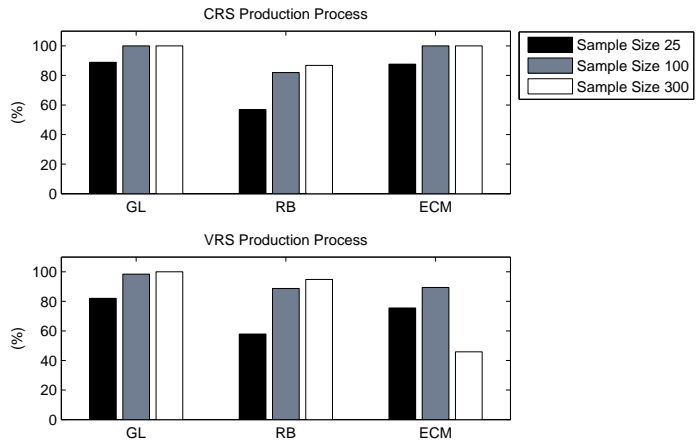}
	\caption{\footnotesize Impact of Variations in Sample Size on the Correct Identification of Efficient DMUs}
	\label{samplesize}
\end{figure}

As different applications involve different sample sizes, we use Experiment 8 and 9 to investigate the impact of small and large sample sizes on the performance of the three variable selection methods.  As can be seen from Table~\ref{tab_perfMetric}, \ref{tab_perfMetric2} and \ref{tab_identification},  the GL method outperforms both ECM method and RB test regardless of the sample size.  The outperformance of the GL method is more evident for a VRS production process.  For instance, while the ECM method and RB test identify 75\% and 58\% of efficient DMUs as efficient respectively, the GL method correctly identifies 85\% of them.  These observations can be summarized in Figure~\ref{samplesize}.  From the figure, it is clear that the RB test gets most affected by the small sample size for both CRS and VRS production processes.  The ECM method exhibits somewhat strong performance except when the sample size is increased for a VRS production process.  In general, consistent with Nataraja and Johnson's \citeyearpar{JohnsonNataraja} results, the performance of the three methods improves as the sample size increases.

On the other hand, an increase in sample size negatively influences the run time of each variable selection method.  In Experiment 9 with a large sample size of 300, while the execution time only takes 17.09 (25.08) seconds for the GL method for a CRS (VRS) production process, it takes 904.45 (314.51) and 114.30 (97.47) seconds for the ECM method and RB test respectively (see Table~\ref{tab_CPU}).  Across all 12 experiments, the GL method has the shortest execution time while the ECM method, which uses the backward selection algorithm, has the longest execution time.
\begin{table}[htbp]
\begin{center}
\caption{\footnotesize CPU Time (seconds)}
\begin{minipage}[t]{25pc}
\scriptsize Table~\ref{tab_CPU} presents the amount of time each method took for one trial of each experiment.
\begin{center}
{\footnotesize\tabcolsep=5pt\begin{tabular*}{25pc}{@{\extracolsep{\fill}}  c | c c c | c c c }
\hline%
\multicolumn{1}{c}{RTS:} & \multicolumn{3}{c} {CRS} &\multicolumn{3}{c}{VRS}\\
\hline%
Experiment & GL & ECM & RB & GL & ECM & RB \\
\hline
1&1.62&27.02&4.14&2.12&27.51&4.14\\
2&1.54&26.97&4.38&2.17&33.23&5.11\\
3&1.39&29.57&2.78&1.85&33.73&3.12\\
4&1.94&24.20&4.69&1.85&30.12&5.35\\
5&1.65&28.30&4.17&2.18&32.88&3.40\\
6&1.51&28.41&3.65&1.86&32.66&4.50\\
7&1.66&24.77&1.49&2.19&27.40&3.69\\
8&0.36&2.42&0.34&0.35&2.68&0.34\\
9&17.09&904.45&114.30&25.08&314.51&97.47\\
10&1.95&33.80&4.53&2.96&44.91&4.47\\
11&1.77&16.83&3.27&1.98&18.19&4.14\\
12&2.58&60.98&4.41&2.84&62.32&4.50\\
\hline
\end{tabular*}}%
\label{tab_CPU}
\end{center}
\end{minipage}
\end{center}
\end{table}%

\subsection{Impact of Variations in the Importance of Inputs in the Production Process}
\begin{figure}[h!]
	\centering
	\includegraphics[height = 2.5in]{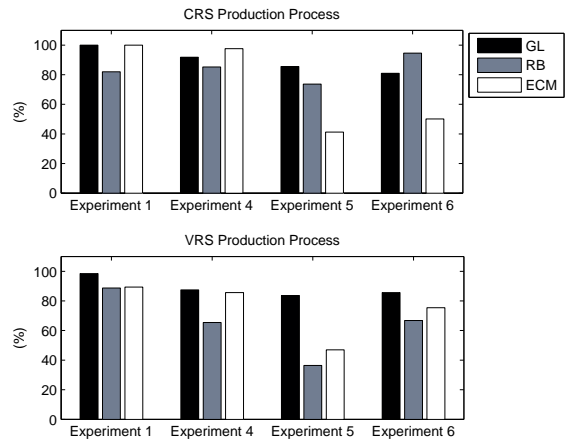}
	\caption{\footnotesize Impact of Variations in Input Contribution on the Correct Identification of Efficient DMUs}
	\label{inputcontribution}
\end{figure}

As it is reasonable to assume considerable variations in the relative importance of inputs in the production process, it is essential to test the robustness of the results with respect to the variations in input contribution to output.  When input contribution to output is varied (i.e Experiments 4, 5, and 6), all three variable selection methods, in general, have better performance under a CRS production process.  As can be seen from Figure~\ref{inputcontribution}, which plots the results of Experiments 1, 4, 5, and 6, the GL method exhibits the most stable and strongest performance across the four experiments.  Although the RB test outperforms the GL method in Experiment 6 under a CRS production process, it underperforms under a VRS production process.  The performance of both ECM method and RB test gets negatively impacted by the variations in input contribution to output when inputs are correlated, especially under a VRS production process.  For example, in terms identifying efficient and inefficient DMUs, when input contribution is varied for correlated inputs (i.e. Experiment 5), the ECM method and RB test identify 36\% and 47\% of efficient DMUs correctly under a VRS production process.  These values are 5\% and 27\% lower than the corresponding values obtained for a CRS production process.  Similar observations can be made in terms of MSEs between the true and model efficiency estimates.   Consistent with the results obtained so far, the GL method shows equal or better performance compared to the other two benchmark methods.

\subsection{Impact of Varitions in the Dimensionality of the Production Process}
Experiments 10 and 11 consider variations in the dimensionality of the production process.  When the dimensionality of the production function is increased in Experiment 10, the GL method outperforms both ECM method and RB test.  For a CRS production process, the RB test has the weakest performance in terms of MSE and correlations between the true and model efficiency estimates.  Under a VRS production process, both ECM method and RB test show considerably weaker performance.  More specifically, MSEs and correlations between the true and model efficiency estimates are at least 10 times higher and 8\% lower than those obtained for the GL method respectively.  As the dimensionality of the production function decreases in Experiment 11, the performance of all three methods improves.

From the results obtained from various experiments in our simulation study, we can conclude that the GL method significantly outperforms both ECM method and RB test and is most robust to variations in the covariance structure of inputs, sample size, importance of inputs, and the dimensionality of the production space.  Also, the GL method is the fastest algorithm and is least vulnerable to the choice of underlying production technology among the three variable selection methods presented.

\section{Conclusion}
\label{sec:conclusion}
As Golany and Roll \citeyearpar{GolanyRoll} note in their study, surprisingly a few number of studies give an overview of DEA as an application procedure that must focus on the choice of variables.  Even though no functional form is specified in DEA models, DEA results heavily rely on the selection of input and output variables.  Wrong choices of variables are likely to compromise the accuracy of the analysis.  Likewise, model specification must be a central concern in DEA.

In this study, we have developed a data-driven joint variable selection method based on group Lasso and reported its significant outperformance over the prevailing variable selection methods, the efficiency contribution measure (ECM) method and regression-based (RB) test.  We have evaluated the performance of our proposed method and its benchmarks by means of Monte Carlo simulations and examined the sensitivity of the results to the variations in the correlations between inputs, sample size, importance of inputs, and the dimensionality of the production space.  Based on the results obtained from a diversified set of simulation experiments, it is evident that the GL method is more robust and efficient than its benchmarks.  To conclude, this study proposes a more sophisticated and quantitative variable selection method that will help finding a parsimonious DEA model, which uses as many variables as needed, but as few as possible.

\bibliographystyle{abbrvnat}
\bibliography{bibliography}

%
%
%





\end{document}

%% file: rax_pdfsetup.tex
\usepackage{hyperref}

\hypersetup{
pdfauthor = {Irene Song},
pdftitle = {\maintitle},
pdfsubject = {\subject},
pdfkeywords = {\keywordlist}}

%% file: tony_def.tex
\newtheorem{thm}{Theorem}[section]

\newcommand{\barV}{\bar{V}}
\newcommand{\barX}{\bar{X}}
\newcommand{\barY}{\bar{Y}}
\newcommand{\barI}{\bar{I}}
\newcommand{\barv}{\bar{v}}

\newcommand{\transp}{^\top}

\newcommand{\calL}{\mathcal{L}}

\newcommand{\calT}{\mathcal{T}}

\newcommand{\setTo}[4]{\left\{#1\right\}_{#2=#3}^#4}

\newcommand{\sumTo}[3]{\sum_{#1=#2}^{#3}}

\newcommand{\eqnBegin}{\begin{equation}}
\newcommand{\eqnEnd}{\end{equation}}
\newcommand{\eqnSBegin}{\begin{equation*}}
\newcommand{\eqnSEnd}{\end{equation*}}
\newcommand{\eqnarrayBegin}{\begin{eqnarray}}
\newcommand{\eqnarrayEnd}{\end{eqnarray}}
\newcommand{\eqnarraySBegin}{\begin{eqnarray*}}
\newcommand{\eqnarraySEnd}{\end{eqnarray*}}

\newcommand{\varpone}[2]{#1^{#2+1}}

\newcommand{\xkpone}{\varpone{x}{k}}
\newcommand{\ykpone}{\varpone{y}{k}}

\newcommand{\xsupk}[2]{#1^{(#2)}}

\newcommand{\gammakpone}{\xsupk{\gamma}{k+1}}
\newcommand{\gammak}{\xsupk{\gamma}{k}}